\documentclass[12pt, a4paper]{article}

\usepackage{tikz}

\usepackage{my_proof}

\usepackage{amsmath}                                                            
\usepackage{amssymb}                                                            
\usepackage{amsfonts}                                                           

\usepackage{epsfig}                                                             
                                                                                
\usepackage{ifthen}                                                             

\newcounter{bsplemsatz}
\renewcommand{\thebsplemsatz}{\arabic{bsplemsatz}}

\newenvironment{env}[2][]{%
  \refstepcounter{bsplemsatz}                                                   
  \trivlist                                                                     
  \item[\hskip\labelsep{\bf #2\relax\ \thebsplemsatz.}%
  ]\ignorespaces\itshape%
}{%
  \endtrivlist                                                                  
}                                                                               

\newcommand{\NN}{\mathbb{N}}                                                    
\newcommand{\ZZ}{\mathbb{Z}}                                                    
                                                    
\newcommand{\RR}{\mathbb{R}}                                                    
\newcommand{\CC}{\mathbb{C}}                                                    

\DeclareMathAlphabet{\masc}{U}{eus}{m}{n}                                       
\DeclareMathAlphabet{\mafr}{U}{euf}{m}{n}                                       
\DeclareMathAlphabet{\malf}{OT1}{cmtt}{m}{it}               

\newcommand{\Ker}{\operatorname*{ker}}
\newcommand{\rk}{\operatorname*{rk}}

\newcommand{\ord}[2]{{\operatorname*{ord}}_{#1}(#2)}

\newcommand{\notdiv}{\!\nmid\!}
\newcommand{\col}{\operatorname*{col}}

\newcommand{\I}{{\mafr i}}

\begin{document}

\title{On the eigenvalues of distance powers of circuits}

\author{%
J. W. Sander\\
Institut f\"ur Mathematik und Angewandte Informatik\\ Universit\"at Hildesheim\\ D-31141 Hildesheim\\ Germany\\
\and 
T. Sander\\
Institut f\"ur Mathematik\\ Technische Universit\"at Clausthal\\ D-38678 Clausthal-Zellerfeld\\ Germany\\
}

\date{\today}

\maketitle



\thispagestyle{empty}

\centerline{\large \bf Abstract}
Taking the $d$-th distance power of a graph, one adds edges between all pairs of vertices
of that graph whose distance is at most $d$. It is shown that only the numbers
$-3$, $-2$, $-1$, $0$, $1$, $2d$ can be
integer eigenvalues of a circuit distance power. Moreover, their respective
multiplicities are determined and explicit constructions for corresponding eigenspace bases
containing only vectors with entries $-1$, $0$, $1$ are given.

{\bf Keywords:} distance power, circuit graph, integral eigenvalues, eigenspace basis

{\bf 2000 Mathematics Subject Classification:} Primary 05C50, Secondary 15A18\\ 




\section{Introduction}

Given a graph $G$ and a positive integer $d$, the $d$-th {\em distance power} $G^{(d)}$ of $G$ is obtained from $G$ by
adding edges between all pairs of vertices whose distance is at most $d$. This implies that
$G^{(1)}$ is isomorphic to $G$.
We are interested in distance powers of the circuit graph on $n$ vertices (denoted by $C_n$).
They belong to the important class of {\em circulant graphs}.

Circulant graphs are characterized as follows. 
Assume that the vertices of a given graph are $0, 1, \ldots, n-1$ and
consider the set $N$ of neighbors of vertex $0$. The graph is circulant
if and only if under every possible cyclic rotation of the vertex numbers the set of
neighbors of the new vertex $0$ remains $N$. We shall call $N$ the {\em jump set} of the graph.
Circulant graphs have many fascinating properties, cf.~{\sc Boesch, Tindell} \cite{boesch84}, 
and interesting applications.  For example, they play a role in the study of
redundant communication networks ({\sc Bermond et al.} \cite{bermond95}).

Moreover, circulant graphs model quantum systems. Such a system
is periodic if and only if its graph is {\em integral}, i.e., if it has only integer eigenvalues
({\sc Saxena et al.} \cite{saxsevsh07}).
The eigenvalues of a graph are the eigenvalues of its adjacency matrix. This matrix is 
defined by numbering the vertices of the graph with $0,\ldots,n-1$ and 
letting the entry at position $(i,j)$ be one if the vertices numbered $i$ and $j$ are adjacent
and zero otherwise\footnote{When dealing with circulant graphs it is convenient to
use zero-based matrix indices.}. The zero-one pattern of the adjacency matrix depends on the chosen vertex numbering 
of the vertices, but it follows from basic linear algebra that its eigenvalues do not.
The set of eigenvalues of a graph, called its {\em spectrum},
reflects several structural properties of the graph (see e.g. the books \cite{eigenspacesofgraphs}, \cite{cvetdoob82}).

There exists an elegant condition due to {\sc So} \cite{so06} that asserts a given circulant graph 
with vertices $0, 1, \ldots, n-1$ is integral if and only if
its jump set $N$ consists of complete sets of numbers having the same gcd with $n$.
Let us check this condition for distance powers $C_n^{(d)}$ of circuits.
Clearly, the jump set of a distance power $C_n^{(d)}$ is $\{1,2,\ldots,d,n-d,\ldots,n-2,n-1\}$.
Note that for odd $n$ the condition $k\in N$ is equivalent to $C_n^{(d)}$ being the complete graph,
for even $n$ take $k-1\in N$.
Now assume $n=2k+1$ and let $g=\gcd(k,n)$. From the properties of the gcd it easily follows that necessarily $g\vert 1$
and therefore $k\in N$, since $d\geq 1$.
Similarly, it follows for $n=2k$ and $g=\gcd(k-1,n)$ that $g\vert 2$. So, for $d\geq 2$, we have $k-1 \in N$.
The bottom line is that, except for a few trivial cases with $d=1$, integral distance powers are necessarily complete graphs.
It is readily checked that the only non-complete circuit distance power is $C_6$, for example by using 
the well known explicit eigenvalue formula for $C_n$ (cf.~{\sc Biggs} \cite{biggs}).

So, since integrality is out of reach for non-complete circuit distance powers, we answer the question which integers are possible eigenvalues at all. Some partial results exist in the literature.
For instance, it is well known and easy to show that the circuit
graph $C_n$ itself has eigenvalue $0$ if and only if $4\vert n$ (see {\sc Biggs} \cite{biggs}).
Much more involved arithmetic expressions are required to describe singularity of
circuit distance powers $C_n^{(d)}$ (see {\sc Sander} \cite{singdist} and Theorem \ref{cndnullity} in section \ref{secintegereigenvals}). In the special case of circuit squares $C_n^{(2)}$ related results
were found by {\sc Davis et al.} in \cite{davisdom00}.
We show that only the numbers $-3$, $-2$, $-1$, $0$, $1$, $2d$ can be
integer eigenvalues of a $d$-th circuit distance power and determine the associated
eigenvalue multiplicities. This is the first goal of the present work, covered in section \ref{secintegereigenvals}.

The second goal is to study the eigenspaces associated with the integral eigenvalues.
We will show in section \ref{seceigenspaces} that it is always
possible to choose {\em simply structured} bases, in the sense that the basis vectors contain only entries
from the set $\{-1, 0, 1\}$. Such bases have been shown to exist for a number of graph classes.
Usually, attention is restricted to the graph kernel, i.e. the 
eigenspace for the eigenvalue $0$. The existing literature features results
on trees and forests ({\sc John, Schild} \cite{johnschild96},
{\sc Sander, Sander} \cite{fox1}, {\sc Akbari et al.} \cite{akbariali06}), line graphs of trees ({\sc Marino et al.}~\cite{marino06}, {\sc Sciriha} \cite{sciriha99}),
unicyclic graphs ({\sc Nath, Sarma} \cite{nathsarma07}, {\sc Sander, Sander} \cite{fox2}), bipartite graphs ({\sc Cvetkovi\'{c}, Gutman} \cite{cvetgut72}), or
cographs ({\sc Sander} \cite{eigcograph}).
There exists analogous research concerning the incidence matrix of a graph, where 
the problem of finding simple kernel bases can be considered as solved (cf.~{\sc Villarreal} \cite{villa95}, {\sc Hazama} \cite{hazama02} or {\sc Akbari et al.} \cite{akbarigha06}).
What makes circuit distance powers interesting is that we can construct simply structured bases
for {\em all} eigenspaces of integer eigenvalues. Such a property is obvious for the complete graphs $K_n$ and,
hence, for all usual product graphs (cf.~\cite{eigenspacesofgraphs}) that can be derived from them, for example 
Sudoku graphs (cf. {\sc Sander} \cite{sudokuintegral}). However, such products are integral and it seems like
the non-complete circuit distance powers are the first known class of non-integral graphs (excepting $C_6$) with this
property. 

Finally, in section \ref{secoutlook} we consider multiplicities of arbitrary general eigenvalues of circuit distance powers.
We show that all eigenvalues of $C_n^{(d)}$ that lie in the interval $(2d, d/3)$ have multiplicity two.
Moreover, we observe only $0$ and $-2$ may be single eigenvalues. We close with an outlook on path distance powers
where the situation is quite the opposite.

\section{Integer eigenvalues and their multiplicities}\label{secintegereigenvals}

There exists an explicit formula for the eigenvalues of a circuit distance power. The key is the
observation that they belong to the class of circulant graphs.
In this section, we will tacitly assume that all considered circuit powers $C_n^{(d)}$ are non-complete,
i.e.~$1\leq d < \frac{n-1}{2}$.

A matrix in which the $i$-th column vector (counting from $i=0$) can be derived from the first column vector
by means of a downward rotation by $i$ entries is called a {\em circulant matrix} \cite{davis}.
Clearly, with respect to some suitable vertex numbering, every circulant graph has a circulant
adjacency matrix.

In the following, let us abbreviate $\omega_n=e^{\frac{2\pi\I}{n}}$.

\begin{env}{Lemma}\cite{biggs} \label{circulantmatrixspectrum}
Let $(a_1,a_2,\ldots,a_n)^T$ be the first column of a real circulant matrix $A$.
Then the eigenvalues of $A$ are exactly
\begin{equation}\label{eqcircuitpowersingular00}
  \lambda_r = \sum\limits_{j=1}^n a_j \omega_{n}^{(j-1)r} , \quad r=0,\ldots,n-1.
\end{equation}
where $\omega_n=e^{\frac{2\pi\I}{n}}$.
\end{env}

\begin{env}{Theorem}\label{circuitpowerspectrum} The eigenvalues of $C_n^{(d)}$
are exactly
\begin{equation}\label{eqcircuitpowersingular}
 \lambda_0=2d,\quad \lambda_r =  \frac{\sin\left((2d+1)\frac{r}{n}\pi\right)}{\sin\frac{r}{n}\pi}-1
\end{equation}
for $r=1,\ldots,n-1$.
\end{env}

\begin{proof}
Use Lemma \ref{circulantmatrixspectrum} and the
following well-known trigonometric identity for the functions $D_q(x)$ of the {\em Dirichlet kernel} \cite{zygmund55}:
\begin{equation}\label{eqcircuitpowersingular01}
   D_q(x) = \sum\limits_{j=-q}^q e^{iqx} = \frac{\sin\left((q+\frac{1}{2})x\right)}{\sin\frac{x}{2}}.
\end{equation}
\end{proof}

Let us now investigate which integer eigenvalues a circuit distance power can have.
Writing the second part of (\ref{eqcircuitpowersingular}) as
\begin{equation}\label{eqcndtrigon}
(\lambda_r + 1) {\sin\frac{r}{n}\pi} - {\sin\left((2d+1)\frac{r}{n}\pi\right)} = 0
\end{equation}
we see that it is a trigonometric Diophantine equation of the form
\begin{equation}\label{eqtwosine}
  A \sin 2\pi a + B \sin 2\pi b  = C
\end{equation}
with rational numbers $A, B, C, a, b$.

{\sc Conway} and {\sc Jones} have outlined how to find the solutions for such equations.
Theorem~7 of their paper \cite{conwayjones76} considers the similar case
\begin{equation}\label{eqfourcosine}
  A \cos 2\pi a + B \cos 2\pi b + C \cos 2\pi c + D \cos 2\pi d = E.
\end{equation}

Adapting their results, we get all nontrivial solutions of equation (\ref{eqtwosine}) as follows:
\begin{env}{Theorem}\label{conwayjones}
Consider at most two distinct rational multiples of $\pi$ lying in the interval $(0,\pi/2)$
for which some rational linear combination of their sines, but of no proper subset, is rational.
The only possible linear combinations, up to multiplication with a rational nonzero factor, are:
\begin{equation}\label{eqtwocosinesol}
    \sin \frac{\pi}{6} = \frac{1}{2}, \qquad \sin \frac{3\pi}{10} - \sin \frac{\pi}{10} = \frac{1}{2}.
\end{equation}
\end{env}

\begin{env}{Proposition} 
The set of integer eigenvalues of a circuit distance power $C_n^{(d)}$ is a subset of 
$\{-3, -2, -1, 0, 1, 2d\}$.
\end{env}

\begin{proof}
Consider an integer solution $\lambda_r$ of equation (\ref{eqcndtrigon}) with $0<r<n$.
For $\vert \lambda_r + 1 \vert \geq 3$,
Theorem \ref{conwayjones} implies that the equation has no solutions 
with distinct rational degree sine arguments in the interval $(0,\pi/2)$. 
But even permitting arbitrary rational degree sine arguments does not help,
so that the equation cannot be  solved.
\end{proof}

It is well known that the degree of regularity of a connected regular graph is an eigenvalue
of multiplicity one (cf.~{\sc Biggs} \cite{biggs}). Therefore $2d$ is always a single eigenvalue of $C_n^{(d)}$.
This also follows from Theorem \ref{circuitpowerspectrum}.

Next, we determine when and with which multiplicity the integers $-3$, $-2$, $-1$, $0$, $1$ occur as eigenvalues of circuit distance powers.

\begin{env}{Theorem}\label{cndmultminusone}
Let $g=\gcd(2d+1,n)$. Then the multiplicity of $-1$ as an eigenvalue of $C_n^{(d)}$ equals $g-1$.
\end{env}

\begin{proof} 
For $\lambda_r=-1$, equation (\ref{eqcndtrigon}) simply becomes
\[
{\sin\left((2d+1)\frac{r}{n}\pi\right)} = 0,
\]
so that, equivalently, we need to find all positive integers $r<n$ such that $(2d+1)r=ln$ for some integer $l$. With the coprime integers $d':=(2d+1)/g$ and $n':=n/g$ the last identity becomes $d'r=n'l$. Hence
$l=d'l'$ for a suitable integer $l'$, and therefore $1\le r = n'l' <n$.
This means $1\le r < g$.
\end{proof}

Let $\ord{p}{n}$ denote the order of the prime divisor $p$ with respect to $n$, i.e.~%
\[
   \ord{p}{n} = \max \{ j\in\NN_0 : p^j \vert n\}.
\]

\begin{env}{Theorem}\cite{singdist}\label{cndnullity}
For given $n,d\in\NN$ let $g:=\gcd(n,d)$ and $h:=\gcd(n,d+1)$. 
Then the multiplicity of $0$ 
as an eigenvalue of $C_n^{(d)}$ is
\[
    \begin{cases}
      g-1      & \text{if}\ \ord{2}{d+1}\geq \ord{2}{n}, \\
      g+h-1    & \text{if}\ \ord{2}{d+1}< \ord{2}{n} \ ~\text{and}~\ 2\notdiv d,\\
      g+h-2    & \text{if}\ 2\vert n\ ~\text{and}~\ 2\vert d.
    \end{cases}
\]
\end{env}

\begin{proof}
For $\lambda_r=0$, equation (\ref{eqcndtrigon}) takes the form
\[
   {\sin\frac{r}{n}\pi} = {\sin\left((2d+1)\frac{r}{n}\pi\right)},
\]
so that we need to determine all integers $0 < r < n$ and $l\in\NN_0$ such that
$dr=ln \ \text{or}\ 2(d+1)r=(2l+1)n$. A detailed proof can be found in \cite{singdist}.
\end{proof}

Let us point out that, since circuit squares $C_n^{(2)}$ are 
4-circulant graphs of type $4C_n(1,2)$, Theorem~5 in \cite{davisdom00} proves our Theorem \ref{cndnullity}
in the special case $d=2$.

Note that the terms $h+g-1$ and $h+g-2$ in Theorem \ref{cndnullity} are fairly interesting.
With the greatest common divisors $g:=\gcd(n,d)$ and $h:=\gcd(n,d+1)$ we see that the terms are essentially sums
of two multiplicative objects ---
a somewhat irritating fact for number theorists. 

In the same manner as Theorem \ref{cndnullity} we can prove the conditions for eigenvalue $-2$.

\begin{env}{Theorem}\label{cndmultminustwo}
For given $n,d\in\NN$ let $g:=\gcd(n,d)$ and $h:=\gcd(n,d+1)$. Then the multiplicity of $-2$ 
as an eigenvalue of $C_n^{(d)}$ with $d>1$ is
\[
    \begin{cases}
      h-1      & \text{if}\ \ord{2}{d} \geq \ord{2}{n}, \\
      g+h-1    & \text{if}\ \ord{2}{d} < \ord{2}{n} \ ~\text{and}~\ 2\vert d,\\
      g+h-2    & \text{if}\ 2\vert n\ ~\text{and}~\ 2\notdiv d.
    \end{cases}
\]
\end{env}

\begin{env}{Theorem} A circuit distance power $C_n^{(d)}$ has
eigenvalue $1$ if and only if $6\vert n$ and $d\equiv 1 \mod 6$. 
In this case, the multiplicity of the eigenvalue equals two.
\end{env}

\begin{proof}
For $(\lambda_r + 1) = 2$ we see from Theorem \ref{conwayjones} that
equation (\ref{eqcndtrigon}) can only have a solution if one of the
arguments is a multiple of $\pi/2$. To be precise, the first sine term
must equal $1/2$ and the second sine term must equal $1$. 
This leads to the two solutions $r/n = \pi/6$ or $r/n = 5\pi/6$ for the
first term (recall that $0<r<n$) and to $(2d+1)r/n = \pi + 2k\pi$ with $k\in\ZZ$ for the second term.
Hence, $6\vert n$ and $d\equiv 1 \mod 6$.
\end{proof}

Analogously, we obtain the following theorem:

\begin{env}{Theorem} A circuit distance power $C_n^{(d)}$  has
eigenvalue $-3$ if and only if  $6\vert n$ and $d\equiv 4 \mod 6$. 
In this case, the multiplicity of the eigenvalue equals two.
\end{env}



\section{Eigenspaces for integer eigenvalues}\label{seceigenspaces}

According to {\sc Davis} \cite{davis}, the column vectors of the matrix
\[
  F^\ast=n^{-\frac{1}{2}}\left( \omega_n^{ij}\right)_{i,j=0,\ldots,n-1}\in\CC^{n\times n},
\]
which is the conjugate transpose of the {\em Fourier matrix} $F\in\CC^{n\times n}$,
constitute a complete and universal set of complex eigenvectors for {\em every} circulant matrix $M$ of order $n$.
Moreover, the $r$-th column of $F^\ast$, denoted by $\col(r)$, yields a complex eigenvector for
eigenvalue $\lambda_r$ of Theorem \ref{circuitpowerspectrum}. 

In the following, we use this fact to establish {\em real} eigenspace bases. Even more, we assert
that for all integer eigenvalues of a circuit distance power $C_n^{(d)}$ there exist associated 
simply structured eigenspace bases.

\begin{env}{Theorem}\label{cndsimplystructured} Every integer eigenvalue of 
$C_n^{(d)}$ admits a simply structured eigenspace basis.
\end{env}

\begin{proof}
{\em Case $\lambda_r=2d$:} It is well-known \cite{biggs} that the all-ones vector forms a corresponding eigenspace basis.

{\em Case $\lambda_r=-3$ or $\lambda_r=1$:} It
it easily verified that in both cases the vectors
\[
  (1,1,0,-1,-1,0,\ldots)^T, (1,0,-1,-1,0,1,\ldots)^T 
\]
form a corresponding simply structured basis.

{\em Case $\lambda_r=0$:} Let $g=\gcd(n,d)$ and $h=\gcd(n,d+1)$. It follows from the proof of Theorem \ref{cndnullity}
that the vectors $u_1,\ldots,u_{g-1}$ with $u_k=\sqrt{n}\cdot \col(kn/g)$ form a basis of a subspace of $\Ker C_n^{(d)}$.
We will show that the vectors $u'_1,\ldots,u'_{g-1}$ with 
\[
     u'_k= \sum\limits_{m=0}^{n/g-1} e_{k+mg} -e_{g+mg}
\]
constitute
an alternative (real) basis of this subspace.

Let $M$ be the matrix with columns $u_1,\ldots,u_{g-1}$. 
Fix some $1\leq\iota\leq g-1$ and let $M'$ be the matrix with columns
$u_1,\ldots,u_{g-1},u'_\iota$. Clearly, $\rk M' \geq \rk M = g-1$.
Actually, we have $\rk M' = g-1$ since the sum of all row vectors of $M'$ vanishes.
To see this, consider the summation of the values in a single column.
We have $u_k=(\omega_n^{0kn/g},\omega_n^{1kn/g},\omega_n^{2kn/g},\ldots,\omega_n^{(g-1)kn/g})^T$
so that its component sum is
\begin{equation}\label{eqgaussper1}
    \sum\limits_{m=0}^{g-1} \omega_n^{mkn/g} = \sum\limits_{m=0}^{g-1} \omega_g^{km}
\end{equation}
and therefore a Gaussian period. Because of $1\leq k\leq g-1$ we have $g\notdiv k$ so that,
according to the theory of Gaussian periods (see {\sc Davenport} \cite{davenport80} or {\sc Nagell} \cite{nagell51}), the component sum in equation (\ref{eqgaussper1}) is zero.
Moreover, the component sum of $u'_\iota$ is zero, too.
Hence it follows that $u'_\iota$ is a linear combination of the vectors $u_k$. 
Since the vectors $u'_1,\ldots,u'_{g-1}$ are obviously linearly independent 
we see that they are a basis for the space spanned by $u_1,\ldots,u_{g-1}$.

In the case that $\ord{2}{d+1}< \ord{2}{n}$, equivalently
$2h\vert n$, the vectors  $v_1,\ldots,v_{h}$ with $v_k=\sqrt{n}\cdot \col(kn/h-n/(2h))$
form a basis of another subspace of $\Ker C_n^{(d)}$.
 
A similar argument as for the vectors $u'_k$ shows that the vectors
$v'_1,\ldots,v'_{h}$ with 
\[
    v'_k=\sum\limits_{m=0}^{n/h-1} (-1)^m e_{k+mh}
\]
constitute a basis of the subspace of $\Ker C_n^{(d)}$ spanned by the vectors $v_1,\ldots,v_{h}$.

Let us consider the cases listed in Theorem \ref{cndnullity}:
\begin{itemize}
\item If $\ord{2}{d+1}\geq \ord{2}{n}$, then $\{u'_1,\ldots,u'_{g-1}\}$ is a basis of $\Ker  C_n^{(d)}$.
\item If $\ord{2}{d+1}< \ord{2}{n}$ and $2\notdiv d$, then $\{u'_1,\ldots,u'_{g-1},v'_1,\ldots,v'_{h}\}$ is a basis of $\Ker  C_n^{(d)}$.
\item If $2\vert n$ and $2\vert d$, then $\{u'_1,\ldots,u'_{g-1},v'_1,\ldots,v'_{h}\}$ can be reduced to a basis of $\Ker  C_n^{(d)}$.
\end{itemize}
All of the above bases are simply structured.

{\em Case $\lambda_r=-2$:} Use Theorem \ref{cndmultminustwo}. This case is analogous to case $\lambda=0$, only with swapped roles of $g$ and $h$.
We have complex subspace basis vectors $u_1,\ldots,u_{g}$ with $u_k=\sqrt{n}\cdot \col(kn/g-n/(2g))$ and can find 
real basis vectors $u'_1,\ldots,u'_{g}$ with 
\[
    u'_k=\sum\limits_{m=0}^{n/g-1} (-1)^m e_{k+mg}
\]
for the same subspace.
Likewise, we have complex vectors $v_1,\ldots,v_{h-1}$ with $v_k=\sqrt{n}\cdot \col(kn/h)$ 
and real vectors $v'_1,\ldots,v'_{h-1}$ with
\[
     v'_k= \sum\limits_{m=0}^{n/h-1} e_{k+mh} -e_{h+mh}.
\]

{\em Case $\lambda_r=-1$:} With the help of Theorem \ref{cndmultminusone}, we can reason as in the first part
of case $\lambda=0$, but with $g=\gcd(n,2d+1)$. The same complex vectors $u_1,\ldots,u_{g-1}$ 
and real vectors $u'_1,\ldots,u'_{g-1}$ are obtained.
\end{proof}

\begin{env}{Example}
\newcommand{\m}{\phantom{-}}
Theorem \ref{cndsimplystructured} asserts that the vectors
\[
\begin{split}
   (\m1, \m0,  -1, \m0, \m1, \m0,  -1, \m0, \hdots,  \m1, \m0,  -1, \m0)^T,\\
   (\m0, \m1, \m0,  -1, \m0, \m1, \m0,  -1, \hdots,  \m0, \m1, \m0,  -1)^T,\\
   (\m1, \m0,  -1, \m1, \m0,  -1, \m1, \m0, -1, \hdots,  \m1, \m0, -1)^T,\\
   (\m0, \m1,  -1, \m0, \m1,  -1, \m0, \m1, -1, \hdots,  \m0, \m1, -1)^T\phantom{,}
\end{split}
\]
form a simply structured eigenspace 
basis of $C_{36}^{(14)}$ for eigenvalue $-2$.
\end{env}

\begin{env}{Remark}
Note that some of the  simply structured bases 
constructed in Theorem \ref{cndsimplystructured} are actually orthogonal.

For $\lambda=2d$ and $\lambda=-1$ this is always the case. For $\lambda=-3$ and $\lambda=1$
one can never obtain a simply structured basis. For $\lambda=0$, we see by
Theorem \ref{cndnullity} that the constructed basis is simply structured if 
${\rm ord}_2(d+1) < {\rm ord}_2(n)$ and $\gcd(n,d)\le 2$. 
Theorem \ref{cndmultminustwo} implies an analogous statement for $\lambda=-2$.
\end{env}

\section{Eigenvalue multiplicities in general}\label{secoutlook}

Let us revisit equations (\ref{eqcircuitpowersingular00}) and (\ref{eqcircuitpowersingular01}) by considering the functions
$f_d: [0,2\pi] \rightarrow \RR$ with
\[
   f_d(\varphi) := \frac{\sin((2d+1)\varphi/2)}{\sin(\varphi/2)}
\]
for $\varphi\in(0,2\pi)$ and the continuous extension $f_d(0)=f_d(2\pi)=2d+1$.
 
Let us point out some obvious properties of $f_d$. We have $f_d(\pi)=\sin((2d+1)\pi/2) \in \{1, -1\}$. Moreover,
$f_d(\varphi)$ is axis symmetric with respect to $\varphi=\pi$.
The zeros of $f_d$ are exactly the integer multiples $kq$ of $q:=2\pi/(2d+1)$ with $k=1,\ldots,2d$.

Since we can find the eigenvalues of $C_n^{(d)}$ as $\lambda_0=2d$ and $\lambda_r = f_d(2\pi r/n)-1\not=2d$ for 
$r=1,\ldots,n-1$, the following fact is obvious:

\begin{env}{Observation}\label{cndoddmult}
Any eigenvalue of $C_n^{(d)}$ of odd multiplicity 
must necessarily be
$2d$, $0$ or $-2$.
\end{env}

\begin{env}{Theorem}\label{cndtwofold}
Every eigenvalue of $C_n^{(d)}$ that is greater than $\frac{d}{\pi}-1$ and less than $2d$ has multiplicity two.
\end{env}

\begin{figure}
\epsfig{file=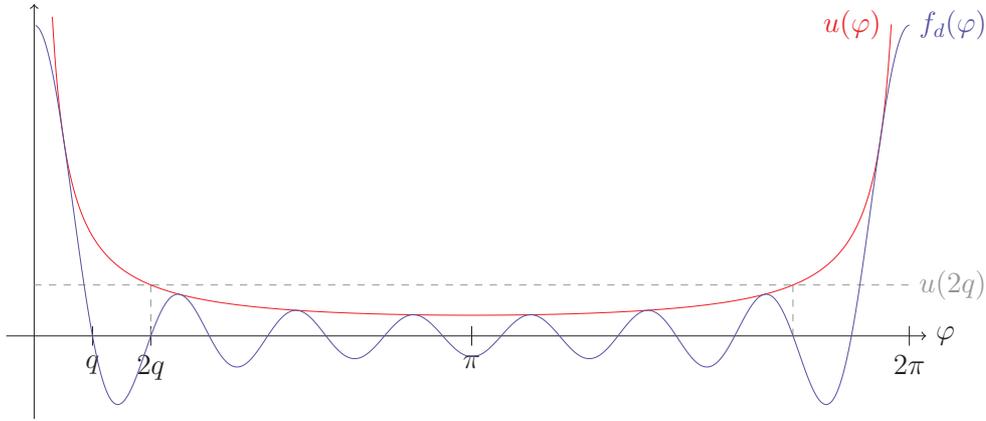,width=14cm}
\caption{Obtaining a lower bound for eigenvalues with multiplicity two }\label{figbound}
\end{figure}

\begin{proof}
A simple upper bound for $f_d$ is 
\[
    u: \varphi \mapsto \frac{1}{\sin(\varphi/2)}.
\]
This bound is strictly decreasing on $(0,\pi)$. 
Observing symmetry, it follows that
$u(2q) \geq f_d(\varphi)$ for $\varphi\in(2q, 2\pi-2q)$.

We have $f_d(0)=2d+1$. So it is clear that $f_d$ is nonnegative in the interval $(0, q)$
and non-positive in the interval $(q, 2q)$.

Claim: $f_d$ is strictly decreasing in the interval $(0, q)$.

We consider the derivative
\[
 f'_d: \varphi \mapsto 
           \frac{(2d+1)\cos\left( (2d+1)\varphi/2 \right) \sin (\varphi/2)%
              - \sin\left( (2d+1)\varphi/2 \right) \cos(\varphi/2)}
           {2\sin^2(\varphi/2)} 
\]
and show that $f'_d(\varphi)\leq 0$ for $0\leq \varphi\leq q$.
It is easy to see that the first term of the numerator is
positive for $0<\varphi<q/2$ and negative for $q/2< \varphi<q$ whereas
the second term of the numerator is positive for $0<\varphi<q$.
Clearly, $f'_d(\varphi)<0$ for $q/2<\varphi<q$. For $0<\varphi<q/2$
consider the ratio of the two numerator terms, which is
\[
    \frac{(2d+1)\cos\left( (2d+1)\varphi/2 \right) \sin (\varphi/2)}%
         {\sin\left( (2d+1)\varphi/2 \right) \cos(\varphi/2)}
  = \frac{ (2d+1)\tan(\varphi/2) }{\tan\left( (2d+1)\varphi/2 \right)}
  \leq 1
\]
as can be concluded, for instance, from the tangent function's Taylor expansion.
This proves the claim.

It follows readily that any eigenvalue of $C_n^{(d)}$
that is greater than $u(2q)-1$ must have multiplicity two, see Figure \ref{figbound}.
An immediate asymptotic relaxation is obtained by observing that $u(2q)>d/\pi$ for all $d\in\NN$ and 
$\lim\limits_{d\rightarrow\infty} \frac{u(2q)}{d/\pi} = 1$.
\end{proof}

If we denote by $\varphi_0$ the unique $\varphi \in (2q,3q)$ such that
$f_d(\varphi)$ is maximal, then it is clear that one could improve the bound $d/\pi-1$ 
in Theorem \ref{cndtwofold} to a bound $d/u(\varphi_0)-1$. To do so, however, requires 
tedious calculations and we could not get anything considerably smaller than $4d/15-1$.

For fixed $1\leq d\leq 4$, the respective smallest graphs with eigenvalues that fulfil the condition of
Theorem \ref{cndtwofold} are $C_5^{(1)}$, $C_7^{(2)}$, $C_{10}^{(3)}$, and $C_{12}^{(4)}$.

We conclude with some remarks on path distance powers $P_n^{(d)}$.

Clearly, one may remove a suitable number of consecutive vertices from a 
circuit distance power to obtain a path graph. It follows from this
observation that there is a certain relation between the spectra of 
path and circuit distance powers since the eigenvalues of a graph
and an induced subgraph interlace (cf. {\sc Haemers} \cite{haemers95}).

However, path distance powers are not circulant and possess many spectral properties
that are quite unlike those of circuit distance powers.
With respect to the questions considered so far in this and the previous sections,
we pose a number of conjectures which we derived from computer experiments:

\begin{env}{Conjecture}
For every integer $k$ there exists a pair $(n, d)$ such that $k$ is an eigenvalue of $P_n^{(d)}$.
\end{env}

\begin{env}{Conjecture}
The complete graph $K_2$ is the only path distance power with eigenvalue $1$.
\end{env}

\begin{env}{Conjecture}\label{simpleevals}
Every eigenvalue $\lambda\not\in \{-2,-1,0\}$ of $P_n^{(d)}$ is simple.
\end{env}

This is a clear contrast to Observation \ref{cndoddmult} and Theorem \ref{cndtwofold}.
Towards proving Conjecture \ref{simpleevals}, it has be shown that for $\frac{n}{2} < d < n-1$ every multiple eigenvalue $\lambda\not= -1$ of $P_n^{(d)}$ has multiplicity two \cite{sanderdiss}. But this result does not apply for values of $d$ outside this range, cf.~the graph $P_{15}^{(6)}$ with triple eigenvalue zero.




\bibliographystyle{amsplain}
\bibliography{cnk_simple_v1}  

\providecommand{\bysame}{\leavevmode\hbox to3em{\hrulefill}\thinspace}
\providecommand{\MR}{\relax\ifhmode\unskip\space\fi MR }
\providecommand{\MRhref}[2]{%
  \href{http://www.ams.org/mathscinet-getitem?mr=#1}{#2}
}
\providecommand{\href}[2]{#2}
\begin{thebibliography}{10}

\bibitem{akbariali06}
S.~Akbari, A.~Alipour, E.~Ghorbani, and G.~Khosrovshahi,
  \emph{$\{-1,0,1\}$-basis for the null space of a forest}, Linear Algebra
  Appl., vol. 414, 2006, pp.~\mbox{506--511}.

\bibitem{akbarigha06}
S.~Akbari, N.~Ghareghani, G.~B. Khosrovshahi, and H.~R. Maimani, \emph{{T}he
  kernels of the incidence matrices of graphs revisited}, Linear Algebra Appl.,
  vol. 414, 2006, pp.~\mbox{617--625}.

\bibitem{bermond95}
J.-C. Bermond, F.~Comellas, and D.F. Hsu, \emph{Distributed loop computer
  networks: a survey}, J. Parallel Distrib. Comput., vol.~24, 1995,
  pp.~\mbox{2--10}.

\bibitem{biggs}
N.~Biggs, \emph{{A}lgebraic graph theory}, second ed., Cambridge Mathematical
  Library, Cambridge University Press, 1993.

\bibitem{boesch84}
F.~Boesch and R.~Tindell, \emph{Circulants and their connectivities}, J. Graph
  Theory, vol.~8, 1984, pp.~\mbox{487--499}.

\bibitem{conwayjones76}
J.~H. Conway and A.~J. Jones, \emph{Trigonometric diophantine equations ({O}n
  vanishing sums of roots of unity)}, Acta Arithmetica, vol.~30, 1976,
  pp.~\mbox{229--240}.

\bibitem{eigenspacesofgraphs}
D.~Cvetkovi{\'c}, P.~Rowlinson, and S.~Simi{\'c}, \emph{{E}igenspaces of
  graphs}, Encyclopedia of Mathematics and its Applications, vol.~66, Cambridge
  University Press, 1997.

\bibitem{cvetdoob82}
D.~M. Cvetkovi{\'c}, M.~Doob, and H.~Sachs, \emph{{S}pectra of graphs. {T}heory
  and application. {S}econd {E}dition.}, VEB Deutscher Verlag der
  Wissenschaften, Berlin, 1982.

\bibitem{cvetgut72}
D.~M. Cvetkovi{\'c} and I.~M. Gutman, \emph{{T}he algebraic multiplicity of the
  number zero in the spectrum of a bipartite graph}, Mat. Vesnik, vol.~9, 1972,
  pp.~\mbox{141--150}.

\bibitem{davenport80}
H.~Davenport, \emph{{M}ultiplicative number theory. 2nd ed. {R}ev. by {H}ugh
  {L}. {M}ontgomery.}, Graduate Texts in Mathematics, Springer, 1980.

\bibitem{davisdom00}
G.~J. Davis, G.~S. Domke, and C.~R.~Garner Jr., \emph{4-{C}irculant {G}raphs},
  Ars Combin., vol.~65, 2000, pp.~\mbox{97--110}.

\bibitem{davis}
P.~J. Davis, \emph{{C}irculant matrices}, John Wiley \&\ Sons, New
  York-Chichester-Brisbane, 1979.

\bibitem{haemers95}
W.~H. Haemers, \emph{Interlacing eigenvalues and graphs}, Linear Algebra Appl.,
  vol. 226/228, 1995, pp.~\mbox{593--616}.

\bibitem{hazama02}
F.~Hazama, \emph{{O}n the kernels of the incidence matrices of graphs},
  Discrete Math., vol. 254, 2002, pp.~\mbox{165--174}.

\bibitem{johnschild96}
P.~E. John and G.~Schild, \emph{{C}alculating the characteristic polynomial and
  the eigenvectors of a tree}, Match, no.~34, 1996, pp.~\mbox{217--237}.

\bibitem{marino06}
M.~C. Marino, I.~Sciriha, S.~K. Simi{\'c}, and D.~V. To{\v{s}}i{\'c},
  \emph{{M}ore about singular line graphs of trees}, Publications de l'Institut
  Math{\'e}matique (Beograd), vol.~79, 2006, pp.~\mbox{70--85}.

\bibitem{nagell51}
T.~Nagell, \emph{{I}ntroduction to number theory}, Wiley, New York, 1951.

\bibitem{nathsarma07}
M.~Nath and B.~K. Sarma, \emph{On the null-spaces of acyclic and unicyclic
  singular graphs}, Linear Algebra Appl., vol. 427, 2007, pp.~\mbox{42--54}.

\bibitem{fox1}
J.W. Sander and T.~Sander, \emph{{O}n {S}imply {S}tructured {B}ases of {T}ree
  {K}ernels}, AKCE J. Graphs. Combin., vol.~2, 2005, pp.~\mbox{45--56}.

\bibitem{sanderdiss}
T.~Sander, \emph{{E}igenspace {S}tructure of {C}ertain {G}raph {C}lasses}, Ph.
  D. Thesis, TU Clausthal, 2004.

\bibitem{singdist}
\bysame, \emph{Singular distance powers of circuits}, Tokyo J. Math., vol.~4,
  2007, pp.~\mbox{491--498}.

\bibitem{eigcograph}
\bysame, \emph{{O}n {C}ertain {E}igenspaces of {C}ographs}, Electron. J.
  Combin., vol.~15, 2008.

\bibitem{sudokuintegral}
\bysame, \emph{{S}udoku graphs are integral}, Electr. J. Combin., vol.~16,
  2009, Research Note N25, 7 pp. (electronic).

\bibitem{fox2}
T.~Sander and J.W. Sander, \emph{{O}n {S}imply {S}tructured {K}ernel {B}ases of
  {U}nicyclic {G}raphs}, AKCE J. Graphs. Combin., vol.~4, 2007,
  pp.~\mbox{61--82}.

\bibitem{saxsevsh07}
N.~Saxena, S.~Severini, and I.E. Shparlinski, \emph{Parameters of integral
  circulant graphs and periodic quantum dynamics}, Int. J. Quantum Inf.,
  vol.~5, 2007, pp.~\mbox{417--430}.

\bibitem{sciriha99}
I.~Sciriha, \emph{{T}he two classes of singular line graphs of trees}, 5th
  Workshop on Combinatorics (Messina, 1999), Rend. Sem. Mat. Messina Ser. II,
  vol.~5, 1999, pp.~\mbox{167--180}.

\bibitem{so06}
W.~So, \emph{Integral circulant graphs}, Discrete Math., vol. 306, 2006,
  pp.~\mbox{153--158}.

\bibitem{villa95}
R.~H. Villarreal, \emph{{R}ees algebras of edge ideals}, Commun. Algebra,
  vol.~39, 1995, pp.~\mbox{3513--3524}.

\bibitem{zygmund55}
A.~Zygmund, \emph{{T}rigonometrical series. 2nd ed.}, Dover Publications, New
  York, 1955.

\end{thebibliography}

\end{document}